\documentclass[11pt]{article}
\usepackage{amsmath, amssymb, amsthm, graphicx}
\usepackage{comment}
\theoremstyle{definition}
\newtheorem{theorem}{Theorem}[section]
\newtheorem{lemma}[theorem]{Lemma}
\newtheorem{proposition}[theorem]{Proposition}

\newtheorem{Remark}[theorem]{Remark}

\newtheorem{Example}[theorem]{Example}
\newenvironment{example}{\begin{Example}\rm}{\end{Example}}

\numberwithin{equation}{section}

\title{The asymptotic behaviour of recurrence coefficients
for orthogonal polynomials with varying exponential weights}

\author{A.B.J. Kuijlaars and P.M.J. Tibboel
\footnote{The authors are supported by
FWO-Flanders project G.0455.04. The first author
is also supported by K.U. Leuven research grant OT/04/21,
by the Belgian Interuniversity Attraction Pole P06/02,
by the European Science Foundation Program MISGAM,
and by a grant from the Ministry of Education and
Science of Spain, project code MTM2005-08648-C02-01.}
\\[10pt]
{\em Department of Mathematics, Katholieke Universiteit Leuven,} \\
{\em Celestijnenlaan 200B, 3001 Leuven, Belgium} \\
arno.kuijlaars@wis.kuleuven.be, pieter.tibboel@wis.kuleuven.be}
\begin{document}
\maketitle
\begin{abstract}
We  consider orthogonal polynomials $\{p_{n,N}(x)\}_{n=0}^{\infty}$
on the real line with respect to a weight $w(x)=e^{-NV(x)}$ and in
particular the asymptotic behaviour of the coefficients $a_{n,N}$
and $b_{n,N}$ in the three term recurrence $x \pi_{n,N}(x) =
\pi_{n+1,N}(x) + b_{n,N} \pi_{n,N}(x) + a_{n,N} \pi_{n-1,N}(x)$.
For one-cut regular $V$ we show, using the Deift-Zhou method of
steepest descent for Riemann-Hilbert problems, that the diagonal
recurrence coefficients $a_{n,n}$ and $b_{n,n}$ have asymptotic
expansions as $n \to \infty$ in powers of $1/n^2$ and powers of
$1/n$, respectively.
\end{abstract}

\section{Introduction}
We consider the asymptotic behavior  of the recurrence
coefficients $a_{n,N}$ and $b_{n,N}$ in the three-term recurrence
relation
\[ x \pi_{n,N}(x) = \pi_{n+1,N}(x) + b_{n,N} \pi_{n,N}(x) +
a_{n,N} \pi_{n-1,N}(x) \]
for orthogonal polynomials with respect to varying exponential weights.
Here $\pi_{n,N}$ is the $n$-th degree monic orthogonal polynomial
with respect to a varying weight
\[ w_N(x) = e^{-NV(x)} \]
where $V$ is real analytic on $\mathbb R$ with
$\lim\limits_{x\to \pm \infty} \frac{V(x)}{\log(1+x^2)} = +\infty$.
Moreover, $V$ is assumed to be one-cut regular, which means
that the equilibrium measure $d\mu_V = \psi_V(x)dx$ associated
with $V$ is supported on one interval $[a,b]$ where it has the form
\begin{equation} \label{psiV}
    \psi_V(x) \, dx = \sqrt{(b-x)(x-a)} h(x) \chi_{[a,b]}(x) \, dx
    \end{equation}
where $h$ is real analytic, strictly positive on $[a,b]$, and in
addition the inequality (\ref{Lagrange 1}) is strict for $x \in
\mathbb R \setminus [a,b]$. See e.g.\ \cite{APS,BI,Deift,EM,ST} for the
definition of the equilibrium measure and for more information on
the one-cut regular case.

Under these assumptions Deift et al.\ \cite{DKMVZ2} proved that
$a_{n,n}$ and $b_{n,n}$ have asymptotic expansions in powers of $1/n$.
Their approach is based on the Deift-Zhou method of steepest descent
applied to the Riemann-Hilbert problem for orthogonal polynomials of
Fokas, Its, and Kitaev \cite{FIK2}. This method was first introduced
in \cite{DZ} and further developed in \cite{DKMVZ1,DKMVZ2,DVZ} and
many papers since then.

The asymptotic result on the recurrence coefficients
was considerably refined by Bleher and Its \cite[Theorem 5.2]{BI} who showed
for polynomial $V$ that
there exists $\varepsilon > 0$ and real analytic functions
$f_{2k}(s)$, $g_{2k}(s)$, $k=0,1, \ldots$, on $[1-\varepsilon, 1+\varepsilon]$
such that the asymptotic expansions
\begin{align} \label{anN-BI}
    a_{n,N} & \sim f_0\left(\frac{n}{N}\right) +
    \sum_{m=1}^{\infty} N^{-2m} f_{2m}\left(\frac{n}{N}\right) \\
    \label{bnN-BI}
    b_{n,N} & \sim g_0\left(\frac{n+1/2}{N} \right) +
    \sum_{m=1}^{\infty} N^{-2m} g_{2m} \left(\frac{n+1/2}{N} \right)
\end{align}
hold uniformly as $n, N \to \infty$ with $1-\varepsilon \leq n/N \leq
1 + \varepsilon$. These $1/N^2$ expansions are used in \cite{BI}
to prove the $1/N^2$ expansion of the free energy
(a.k.a.\ logarithm of the partition function or Hankel determinant)
of the associated random matrix ensemble in the one-cut regular case,
see also \cite{EM}.

The proof of (\ref{anN-BI}) and (\ref{bnN-BI}) in \cite{BI} is
based on the Deift et al.\ result referred to above, in
combination with so-called string equations. It is of some interest
to find a proof that is based on the Riemann-Hilbert
steepest descent analysis only. Here we do this for the
diagonal case $n = N$, and we obtain the following.

\begin{theorem} \label{Main Theorem}
Let $V$ be real analytic and one-cut regular.
Then there exist constants $\alpha_{2m}$ and $\beta_m$, $m = 1,2, \ldots$ (depending
on $V$) such that $a_{n,n}$ and $b_{n,n}$ have the following asymptotic
expansions as $n \to \infty$:
\begin{align} \label{annbnn}
    a_{n,n} \sim\frac{(b-a)^{2}}{16}+\sum\limits_{m=1}^{\infty} \frac{\alpha_{2m}}{n^{2m}},
    \qquad
    b_{n,n} \sim \frac{b+a}{2} + \sum\limits_{m=1}^{\infty} \frac{\beta_m}{n^m},
\end{align}
where $a$ and $b$ are the endpoints of the support of $\psi_V$.
The first coefficient $\beta_1$ in the expansion for $b_{n,n}$ is given explicitly by
\begin{align} \label{beta1}
    \beta_1 = \frac{1}{2\pi(b-a)}  \left(\frac{1}{h(b)}-\frac{1}{h(a)} \right)
\end{align}
where $h$ is the function appearing in the expression (\ref{psiV})
for the equilibrium measure $\psi_V$ associated with $V$.
\end{theorem}

In our proof of Theorem \ref{Main Theorem} we follow the main lines of the
steepest descent analysis of \cite{DKMVZ2}. We will deduce
that the odd powers in the expansion of $a_{n,n}$ vanish from
the structure of the local Airy parametrices around the endpoints.
The expression (\ref{beta1}) for $\beta_1$ is new, although it is
likely that it can be deduced from the approach of \cite{BI} as well.
The explicit formula (\ref{beta1}) shows that $\beta_1 = 0$ if
and only if $h(a) = h(b)$. It is very easy to construct examples of
one-cut regular $V$ such that $h(a) \neq h(b)$ and so $\beta_1 \neq 0$.
We have thus corrected an error in a paper of Albeverio, Pastur,
and Shcherbina \cite[Theorem 1, formula (1.34)]{APS} who claim that
$\beta_1 = 0$ always in the one-cut regular case.

\begin{example}
We may explicitly check Theorem \ref{Main Theorem} using  Jacobi polynomials
with varying parameters $\alpha = AN$, $\beta= BN$, $A,B > 0$. These polynomials
are orthogonal with weight $(1-x)^{AN}(1+x)^{BN}$ on $[-1,1]$.
The equilibrium measure takes the form (\ref{psiV}) with
\begin{equation} \label{abJacobi}
    a,b = \frac{B^2-A^2 \pm 4\sqrt{(1+A+B)(1+A)(1+B)}}{(2+A+B)^2}
    \end{equation}
and
\begin{equation} \label{hxJacobi}
    h(x) = \frac{2+A+B}{2\pi(1-x^2)},
\end{equation}
see \cite{SUV,KVA}. We are in the one-cut regular case, but for
weights restricted to $[-1,1]$. An analysis of the proof of
Theorem \ref{Main Theorem}, however, will show that the results
(\ref{annbnn})-(\ref{beta1}) remain valid in this case as well.

From the explicit form of the recurrence coefficients for Jacobi
polynomials, see e.g. \cite{Chi,KVA},
\begin{align*}
a_{n,n} & =\frac{4(1+A+B)(1+A)(1+B)}{((2+A+B)^2-\frac{1}{n^{2}})(2+A+B)^{2}} \\
b_{n,n} &=\frac{B^{2}-A^{2}}{(2+A+B)(2+A+B+\frac{2}{n})},
\end{align*}
it is easy to see that (\ref{annbnn}) holds. Using
(\ref{abJacobi})-(\ref{hxJacobi}) we can also ascertain the validity
of (\ref{beta1}).
\end{example}

\section{The Riemann-Hilbert Problem}\label{The Riemann-Hilbert
Problem}
The Riemann-Hilbert problem for orthogonal polynomials was
introduced by Fokas, Its, and Kitaev \cite{FIK2}. It asks for a
$2\times2$ matrix valued function $Y(z)$ satisfying
\begin{align}\label{Y(z)}
\left\{\begin{array}{l}
   Y(z)\textrm{ is analytic in }\mathbb{C}\setminus\mathbb{R}
   \\
   Y_{+}(x)=Y_{-}(x)\begin{pmatrix}
            1 & e^{-NV(x)} \\
            0 & 1 \\
            \end{pmatrix} \textrm{ for }x\in\mathbb{R} \\
   Y(z)=\left(\mathnormal{I}+\mathcal{O}\left(\frac{1}{z}\right)\right)\begin{pmatrix}
   z^{n} & 0 \\
   0 & z^{-n} \\
   \end{pmatrix} \textrm{ as } z\to\infty.
\end{array}\right.
\end{align}

The unique solution of (\ref{Y(z)}) is (see e.g.\ \cite{Deift})
\begin{align}\label{SolutionRH}
Y(z)=\begin{pmatrix} \kappa_{n,N}^{-1}p_{n,N}(z) & \frac{1}{2\pi
i\kappa_{n,N}}\int\limits_{\mathbb{R}}\frac{p_{n,N}(t)}{t-z}dt \\
-2\pi i\kappa_{n-1,N}p_{n-1,N}(z) &
-\kappa_{n-1,N}\int\limits_{\mathbb{R}}\frac{p_{n-1,N}(t)}{t-z}dt
\end{pmatrix}
\end{align}
where $p_{n,N}(x) = \kappa_{n,N} \pi_{n,N}(x)$ is the $n$th degree
orthonormal polynomial.
The recurrence coefficients are expressed as follows in terms
of the solution of the Riemann-Hilbert problem (\ref{Y(z)}), see \cite{Deift,DK}.
\begin{proposition}\label{recurrence proposition}
Let
\begin{align}
Y(z)=\left(\mathnormal{I}+\frac{1}{z}Y_{1}+\frac{1}{z^{2}}Y_{2}+\mathcal{O}\left(\frac{1}{z^{3}}\right)\right)\begin{pmatrix}
   z^{n} & 0 \\
   0 & z^{-n} \\
   \end{pmatrix}
\end{align}
Then
\begin{align} \label{anNinY}
a_{n,N}=\left(Y_{1}\right)_{12}\left(Y_{1}\right)_{21}
\end{align}
and
\begin{align} \label{bnNinY}
b_{n,N}=\frac{\left(Y_{2}\right)_{12}}{\left(Y_{1}\right)_{12}}-\left(Y_{1}\right)_{22}
\end{align}
\end{proposition}

For the remainder of this paper we will take $N=n$. We closely
follow \cite{Deift,DKMVZ2} in applying the Deift-Zhou method of
steepest descent for Riemann-Hilbert problems to (\ref{Y(z)}).

\section{The Deift-Zhou method of steepest descent}
The goal of the Deift-Zhou
method of steepest descent for Riemann-Hilbert problems is to change
the original problem into a problem for which the asymptotics for
$z\to\infty$ are normalised and for which all matrices, jump
matrices and solutions alike, are asymptotically close to the
identity matrix for large $n$ which can be solved iteratively. The
specific details and steps needed to achieve this goal shall be
explained below.

\subsection{The First Step: Transformation $Y\mapsto T$}
The key aspect of the first step of the analysis is the equilibrium
measure $\mu_V$ corresponding to $V$. This equilibrium measure $\mu_V$
is the unique probability measure that satisfies for some $\mathnormal{l}$,
\begin{align}\label{Lagrange 1}
   2\int\log|x-y|^{-1} \, d\mu_{V}(y)+V(x) \geq \mathnormal{l}, &
   \qquad\textrm{for all } x\in\mathbb{R}, \\
   2\int\log|x-y|^{-1}\, d\mu_{V}(y)+V(x)=\mathnormal{l}, &
   \qquad\textrm{for all } x\in \textrm{supp } \mu_V. \label{Lagrange 2}
\end{align}
For the one-cut regular case that we are considering we have that the
support is one interval $[a,b]$ and $d\mu_V(x) = \psi_V(x)\, dx$
as in (\ref{psiV}). In addition the inequality (\ref{Lagrange 1}) is
strict for $x \in \mathbb R \setminus [a,b]$.

Define
\begin{align}
g(z)=\int\log(z-s)\, d\mu_{V}(s)=\int\log(z-s)\psi_V(s)\, ds
\end{align}
and
\begin{align}
\phi(z) & = \pi \int_{b}^{z} \left((s-b)(s-a)\right)^{\frac{1}{2}} h(s) \, ds,
\quad z \in \mathbb C \setminus (-\infty,b] \\
\tilde\phi(z) & = \pi \int_{a}^{z}\left((s-b)(s-a))\right)^{\frac{1}{2}} h(s) \, ds,
\quad z \in \mathbb C \setminus [a,+\infty).
\end{align}
If we now put
\begin{align}
T(z)= e^{-n (\mathnormal{l}/2) \sigma_3} Y(z)
 e^{-ng(z)\sigma_{3}} e^{n (\mathnormal{l}/2) \sigma_3},
\end{align}
where $\sigma_{3} = \left(\begin{smallmatrix} 1 & 0 \\ 0 & -1 \end{smallmatrix}\right)$
is the third Pauli matrix, then $T$ satisfies the
Riemann-Hilbert problem
\begin{align}
\left\{\begin{array}{l} T(z)\textrm{ is analytic in
$\mathbb{C}\setminus\mathbb{R}$}, \\[10pt]
T_{+}(x)=T_{-}(x)J_{T}(x)\textrm{ for $x\in\mathbb{R}$},
\\[10pt]
T(z)=I+\mathcal{O}\left(\frac{1}{z}\right)\textrm{ as
$z\to\infty$},
\end{array}\right.
\end{align}
where
\begin{align} \label{JTdef}
J_{T}(x)=\left\{\begin{array}{l}\begin{pmatrix} 1 & e^{-2n\tilde\phi(x)} \\
0 & 1
\end{pmatrix}\textrm{ for $x<a$}, \\
\begin{pmatrix} e^{2n\phi_{+}(x)} & 1 \\
0 & e^{2n\phi_{-}(x)}
\end{pmatrix}\textrm{ for $x\in (a,b)$}, \\
\begin{pmatrix} 1 & e^{-2n\phi(x)} \\
0 & 1
\end{pmatrix}\textrm{ for $x>b$.}
\end{array}\right.
\end{align}

\begin{comment}
\begin{figure}[!hbp]
\begin{center}
\begin{picture}(200,80)(-0,-30)
   \linethickness{0.4 mm}
   \put(0,0){\line(270,0){270}}
   \put(45,-2.5){$\blacktriangleright$}
   \put(145,-2.5){$\blacktriangleright$}
   \put(225,-2.5){$\blacktriangleright$}
   \put(90,0){\circle*{3}}
   \put(180,0){\circle*{3}}
   \put(90,-10){$a$}
   \put(180,-10){$b$}
   \put(10,20){$\left(\begin{array}{cc} 1 & e^{-2n\tilde\phi} \\ 0 & 1\end{array}\right)$}
   \put(190,20){$\left(\begin{array}{cc} 1 & e^{-2n\phi} \\ 0 & 1\end{array}\right)$}
   \put(90,-30){$\left(\begin{array}{cc} e^{2n\phi_{+}} & 1 \\ 0 & e^{2n\phi_{-}}\end{array}\right)$}
   \end{picture}
   \end{center}
   \caption{The jump matrices $J_T$ for the Riemann-Hilbert problem of $T$.}
   \end{figure}
\end{comment}
Since the inequality in (\ref{Lagrange 1})
is strict for $x < a$ and $x > b$ we have that $\tilde{\phi}(x) > 0$
for $x < a$ and $\phi(x) > 0$ for $x > b$. Thus the jump matrices
for $T$ on $(-\infty, a)$ and $(b,\infty)$ tend to the identity
matrix as $n\to\infty$.

\subsection{The Second Step: Transformation $T\mapsto S$}

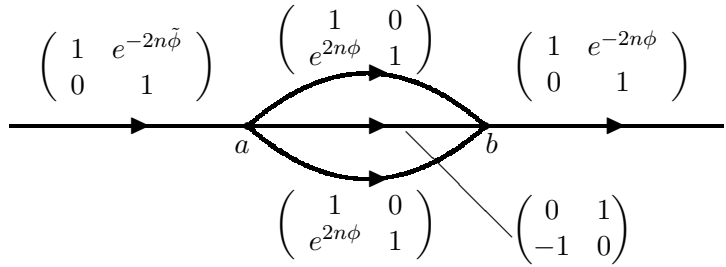
\begin{figure}[!hbp]
\begin{center}
\begin{picture}(270,80)(0,-30)  %(200,80)(0,-30)
   \linethickness{0.4 mm}
   \qbezier(90,0)(135,40)(180,0)
   \qbezier(90,0)(135,-40)(180,0)
   \put(45,-2.5){$\blacktriangleright$}
   \put(135,-2.5){$\blacktriangleright$}
   \put(225,-2.5){$\blacktriangleright$}
   \put(135,17.5){$\blacktriangleright$}
   \put(135,-22.5){$\blacktriangleright$}
   \put(0,0){\line(270,0){270}}
   \put(90,0){\circle*{3}}
   \put(180,0){\circle*{3}}
   \put(85,-10){$a$}
   \put(180,-10){$b$}
   \put(190,-42){$\begin{pmatrix} 0& 1 \\ -1 & 0\end{pmatrix}$}
   \put(10,20){$\left(\begin{array}{cc} 1 & e^{-2n\tilde\phi} \\ 0 & 1\end{array}\right)$}
   \put(190,20){$\left(\begin{array}{cc} 1 & e^{-2n\phi} \\ 0 & 1\end{array}\right)$}
   \put(100,30){$\left(\begin{array}{cc} 1 & 0  \\ e^{2n\phi} & 1 \end{array}\right)$}
   \put(100,-40){$\left(\begin{array}{cc} 1 & 0  \\ e^{2n\phi} & 1 \end{array}\right)$}
   \put(190,-42){\line(-1,1){40}}
   \end{picture}
   \end{center}
   \caption{Jump matrices for $S$ after opening of the lens}\label{TtoS1}
   \end{figure}

The second transformation is the so-called
\textit{opening of the lens} and it is based on the factorisation
\begin{align}\label{factorisation}
\begin{pmatrix} e^{2n\phi_{+}(x)} & 1 \\
0 & e^{2n\phi_{-}(x)}
\end{pmatrix}
=\begin{pmatrix} 1 & 0 \\
e^{2n\phi_{-}(x)} & 1
\end{pmatrix}
\begin{pmatrix}
0 & 1 \\
-1 & 0
\end{pmatrix}
\begin{pmatrix} 1 & 0 \\
e^{2n\phi_{+}(x)} & 1
\end{pmatrix}
\end{align}
of the jump matrix $J_T$ on $(a,b)$. The factorisation
(\ref{factorisation}) allows us to
split the jump on $(a,b)$ as shown in Figure~\ref{TtoS1}.

We use $\Sigma_1$ and $\Sigma_2$ to denote the upper and lower lips
of the lens, respectively. We define $S$ as follows:
\begin{itemize}
\item For $z$ outside the lens, we put $S=T$.
\item For $z$ within the region enclosed by $\Sigma_{1}$ and $(a,b)$,
\begin{align}
S=T\begin{pmatrix} 1 & 0 \\
-e^{2n\phi} & 1
\end{pmatrix}.
\end{align}
\item For $z$ within the region enclosed by $\Sigma_{1}$ and $(a,b)$,
\begin{align}
S=T\begin{pmatrix} 1 & 0 \\
e^{2n\phi} & 1
\end{pmatrix}.
\end{align}
\end{itemize}
Then $S$ satisfies the following Riemann-Hilbert problem:
\begin{align}
\left\{\begin{array}{l} S(z)\textrm{ is analytic in
$\mathbb{C}\setminus(\mathbb{R}\cup\Sigma_{1}\cup\Sigma_{2})$} \\[10pt]
S_{+}(z)=S_{-}(z)J_{S}(z)\textrm{ for
$z\in\mathbb{R}\cup\Sigma_{1}\cup\Sigma_{2}$} \\[10pt]
S(z)=I+\mathcal{O}\left(\frac{1}{z}\right)\textrm{ for $z\to\infty$}
\end{array}\right.
\end{align}
where
\begin{align}
J_{S}(z)=\left\{\begin{array}{l}
\begin{pmatrix} 1 & 0 \\
e^{2n\phi(z)} & 1
\end{pmatrix}\textrm{ for $z\in\Sigma_{1}\cup\Sigma_{2}$},  \\
\begin{pmatrix} 0 & 1 \\
-1 & 0
\end{pmatrix}\textrm{ for $z\in (a,b)$}, \\
\begin{pmatrix} 1 & e^{-2n\tilde\phi(z)} \\
0 & 1
\end{pmatrix}\textrm{ for $z<a$}, \\
\begin{pmatrix} 1 & e^{-2n\phi(z)} \\
0 & 1
\end{pmatrix}\textrm{ for $z>b$},
\end{array}\right.
\end{align}

We may (and do) assume that the lips of the
lens are in the region where $\textrm{Re } \phi < 0$, so that the jump
matrices for $S$ on $\Sigma_{1}$ and $\Sigma_2$ tend to the
identity matrix as $n \to \infty$.

\subsection{The Third Step: Parametrix Away From Endpoints}
The parametrix away from the branch points is  a 'global
solution' $N(z)$ satisfying the Riemann-Hilbert problem
\begin{align}
\left\{\begin{array}{l} N(z)\textrm{ is analytic in
$\mathbb{C}\setminus [a,b]$} \\
N_{+}(x)=N_{-}(x)\begin{pmatrix} 0 & 1 \\
-1 & 0 \end{pmatrix}\textrm{ for $x\in (a,b)$} \\
N(z)=I+\mathcal{O}\left(\frac{1}{z}\right)\textrm{ for
$z\to\infty$}
\end{array}\right.
\end{align}
which has solution (see \cite{Deift})
\begin{align} \label{Ndef}
N(z)&=\begin{pmatrix}\frac{\beta(z)+\beta^{-1}(z)}{2} & \frac{\beta(z)-\beta^{-1}(z)}{2i} \\
   -\frac{\beta(z)-\beta^{-1}(z)}{2i} & \frac{\beta(z)+\beta^{-1}(z)}{2}\end{pmatrix}
\end{align}
where $\beta(z)=\left(\frac{z-b}{z-a}\right)^{\frac{1}{4}}$.

\subsection{The Fourth Step: Parametrices Near Endpoints}

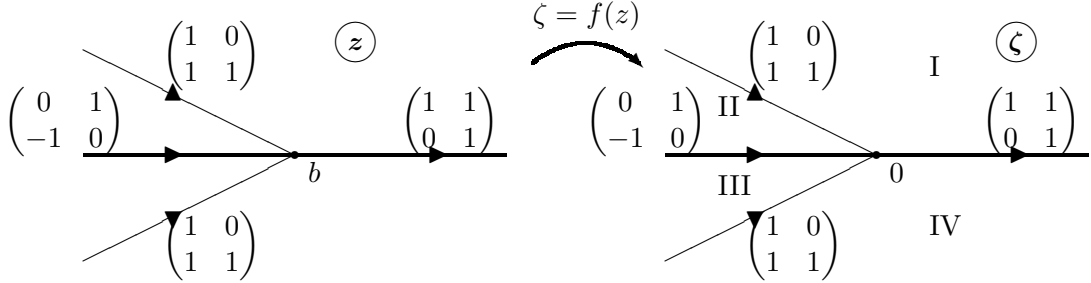
\begin{figure}[t] %[!hbp]
\begin{center}
\begin{picture}(200,80)(0,-30)
   \linethickness{0.4 mm}
   \put(-10,0){\line(270,0){80}}
   \put(-10,0){\line(-2,1){80}}
   \put(-90,0){\line(150,0){80}}
   \put(-10,0){\line(-2,-1){80}}
   %\put(100,25){\vector(2,-1){5}}
   \put(-60,21){$\blacktriangle$}
   \put(-60,-2.5){$\blacktriangleright$}
   \put(40,-2.5){$\blacktriangleright$}
   \put(-60,-27.5){$\blacktriangledown$}
   \put(-5,-10){$b$}
   \put(-120,10){$\begin{pmatrix} 0 & 1 \\
   -1 & 0
   \end{pmatrix}$}
   \put(30,10){$\begin{pmatrix} 1 & 1 \\
   0 & 1
   \end{pmatrix}$}
   \put(-60,35){$\begin{pmatrix} 1 & 0 \\
   1 & 1
   \end{pmatrix}$}
   \put(-60,-37){$\begin{pmatrix} 1 & 0 \\
   1 & 1
   \end{pmatrix}$}
   \put(-10,0){\circle*{3}}
   \put(210,0){\line(270,0){80}}
   \put(210,0){\line(-2,1){80}}
   \put(130,0){\line(150,0){80}}
   \put(210,0){\line(-2,-1){80}}
   %\put(100,25){\vector(2,-1){5}}
   \put(160,21){$\blacktriangle$}
   \put(160,-2.5){$\blacktriangleright$}
   \put(260,-2.5){$\blacktriangleright$}
   \put(160,-27.5){$\blacktriangledown$}
   \put(215,-10){$0$}
   \put(100,10){$\begin{pmatrix} 0 & 1 \\
   -1 & 0
   \end{pmatrix}$}
   \put(230,30){I}
   \put(250,10){$\begin{pmatrix} 1 & 1 \\
   0 & 1
   \end{pmatrix}$}
   \put(150,15){II}
   \put(160,35){$\begin{pmatrix} 1 & 0 \\
   1 & 1
   \end{pmatrix}$}
   \put(150,-15){III}
   \put(160,-37){$\begin{pmatrix} 1 & 0 \\
   1 & 1
   \end{pmatrix}$}
   \put(230,-30){IV}
   \put(210,0){\circle*{3}}
   \qbezier(80,35)(100,50)(120,35)
   \put(120,35){\thicklines\vector(1,-1){2}}
   \put(80,50){$\zeta = f(z)$}
   \put(10,40){\small $\boldsymbol z$}
   \put(13,43){\circle{15}}
   \put(260,40){\small $\boldsymbol \zeta$}
   \put(263,43){\circle{15}}
   \end{picture}
   \end{center}
   \caption{Mapping of neighbourhood of $b$ onto a neighbourhood of $f(b)=0$}
   \label{Conform}
\end{figure}

Having constructed the 'global solution', the next step is finding
'local solutions' close to the endpoints $a$ and $b$. Near $b$, the
local situation is described as in the left picture of
Figure~\ref{Conform} with jump matrix
\begin{align*}
J_{P}(z)= J_S(z) = \left\{\begin{array}{l}
\begin{pmatrix}
1 & 0 \\
e^{2n\phi(z)} & 1
\end{pmatrix}\textrm{ on $\Sigma_{1}\cap U$ and $\Sigma_{2}\cap U$} \\
\begin{pmatrix}
0 & 1 \\
-1 & 0
\end{pmatrix}\textrm{ on  $(a,b)\cap U$} \\
\begin{pmatrix}
1 & e^{-2n\phi(z)} \\
0 & 1
\end{pmatrix}\textrm{ on $(b, \infty)\cap U$}
\end{array}\right.
\end{align*}
where $U$ is a (small) disk around $b$.

\begin{comment}
\begin{figure}[!hbp]
\begin{center}
\begin{picture}(200,80)(0,-30)
   \linethickness{0.4 mm}
   \put(150,0){\line(270,0){100}}
   \put(150,0){\line(-2,1){100}}
   \put(50,0){\line(150,0){100}}
   \put(150,0){\line(-2,-1){100}}
   %\put(100,25){\vector(2,-1){5}}
   \put(100,21){$\blacktriangle$}
   \put(100,-2.5){$\blacktriangleright$}
   \put(200,-2.5){$\blacktriangleright$}
   \put(100,-27.5){$\blacktriangledown$}
   \put(155,-10){$b$}
   \put(100,28){$\mathnormal{\Sigma}_{1}$}
   \put(100,-37){$\mathnormal{\Sigma}_{2}$}
   \put(150,0){\circle*{3}}
   \end{picture}
   \end{center}
   \caption{Sketch of curves on which we have jumps close to
   $b$}\label{localP}
\end{figure}
\end{comment}

We therefore want to find a matrix
function $P$, that solves
\begin{align*}
\left\{\begin{array}{l}
P(z) \textrm{ is analytic on }
U\setminus\left(\Sigma_{1}\cup\Sigma_{2}\cup\mathbb R\right) \\
P_{+}(z)=P_{-}(z)J_{P}(z)\textrm{ on }
(\Sigma_{1}\cup\Sigma_{2}\cup \mathbb R) \cap U \\
P(z)=N(z)\left(\mathnormal{I}+\mathcal{O}\left(\frac{1}{n}\right)\right)\textrm{
as }n\to\infty\textrm{ uniformly for $z\in \partial U$}
\end{array}\right.
\end{align*}
Then $P(z)e^{n\phi(z)\sigma_{3}}$ should have constant jumps on
$(\Sigma_{1}\cup\Sigma_{2}\cup \mathbb R) \cap U$, namely
\begin{align*}
\left(P(z)e^{n\phi(z)\sigma_{3}}\right)_{+}=
\left(P(z)e^{n\phi(z)\sigma_{3}}\right)_{-} \times \left\{\begin{array}{l}
\begin{pmatrix} 1 & 0 \\
1 & 1
\end{pmatrix}\textrm{ for $z\in\left(\Sigma_{1}\cup\Sigma_{2}\right)\cap U$} \\
\begin{pmatrix} 0 & 1 \\
-1 & 0
\end{pmatrix}\textrm{ for $z\in(a,b)\cap U$} \\
\begin{pmatrix} 1 & 1 \\
0 & 1
\end{pmatrix}\textrm{ for $z\in(b,\infty)\cap U$}
\end{array}\right.
\end{align*}

Shrinking $U$ if necessary, we have that
\[ \zeta = f(z) = \left(\frac{3}{2} \phi(z)\right)^{2/3} \]
defines a conformal map from $U$ to a convex neighborhood
of $\zeta = 0$. We may and do assume that the lips of the lens
are taken so that $\Sigma_1 \cap U$ is mapped into $\arg \zeta = 2 \pi/3$,
and $\Sigma_2 \cap U$ is mapped into $\arg  \zeta = 2 \pi/3$,
see Figure~\ref{Conform}.
Denoting the sectors in the $\zeta$-plane by I, II, III, IV
as in Figure~\ref{Conform}, and using the usual Airy function $\textrm{Ai}(\zeta)$,
we construct the Airy model solution $\Phi$ by
\begin{align*}
\Phi(\zeta)&= \left\{ \begin{array}{cl}
\begin{pmatrix}
\textrm{Ai}(\zeta) & \omega \textrm{Ai}(\omega\zeta) \\
\textrm{Ai}'(\zeta) & \omega^{2}
\textrm{Ai}'(\omega\zeta)\end{pmatrix}
    & \textrm{for $\zeta$ in sector IV} \\
\begin{pmatrix}  \textrm{Ai}(\zeta) &-\omega^{2}\textrm{Ai}(\omega^{2}\zeta) \\
\textrm{Ai}'(\zeta) &-\omega \textrm{Ai}'(\omega^{2}\zeta)
\end{pmatrix}
    & \textrm{for $\zeta$ in sector I} \\
\begin{pmatrix}-\omega \textrm{Ai}(\omega\zeta) & -\omega^{2}\textrm{Ai}(\omega^{2}\zeta) \\
-\omega^{2} \textrm{Ai}'(\omega\zeta) & -\omega
\textrm{Ai}'(\omega^{2}\zeta)
\end{pmatrix}
    & \textrm{for $\zeta$ in sector II} \\
\begin{pmatrix}
-\omega^{2}\textrm{Ai}(\omega^{2}\zeta) & \omega
\textrm{Ai}(\omega\zeta) \\
-\omega \textrm{Ai}'(\omega^{2}\zeta) & \omega^{2}
\textrm{Ai}'(\omega\zeta)
\end{pmatrix}
    & \textrm{for $\zeta$ in sector III}
    \end{array} \right. \end{align*}
which has the jump matrices in the $\zeta$-plane
indicated in the right side of Figure~\ref{Conform}.

Then for any analytic prefactor $E_n(z)$ we have that
\begin{align}\label{Phat}
P(z)=E_{n}(z)\Phi(n^{\frac{2}{3}} f(z))e^{n\phi(z)\sigma_{3}}
\end{align}
has the required jump matrices $J_P$. If we choose
\begin{align} \label{En}
E_{n}=\sqrt{\pi}N(z)\begin{pmatrix} 1 & -1 \\
-i & -i
\end{pmatrix} \left(n^{2/3} f(z)\right)^{\sigma_{3}/4}
\end{align}
then the matching condition $P(z) = N(z)(I + \mathcal{O}(1/n))$ as $n \to \infty$
for $z \in \partial U$, is satisfied as well, see e.g.\
\cite{BK,Deift,DKMVZ2} for further detail.

A similar construction yields a parametrix $\tilde P$ in a small
disc $\tilde U$ around $a$. One can see that $\tilde P$ can be obtained
by taking $P$ and interchanging $a$ and $b$ and conjugating with
$\sigma_{3}$.

\subsection{The Fifth Step: Transformation $S\mapsto R$}
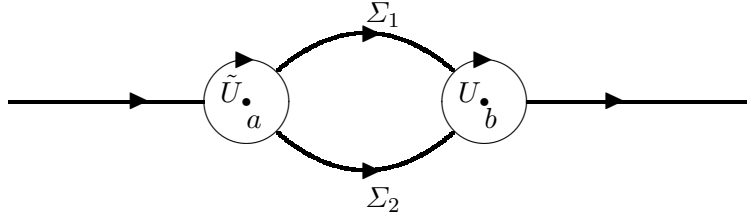
\begin{figure}[!htbp]
\begin{center}
\begin{picture}(280,80)(0,-30)  %(200,80)(0,-30)
   \linethickness{0.4 mm}
   %\qbezier(90,0)(135,40)(180,0)
   \qbezier(102,12)(135,40)(168,12)
   %\qbezier(90,0)(135,-40)(180,0)
   \qbezier(102,-12)(135,-40)(168,-12)
   \put(0,0){\line(270,0){74}}
   %\put(106,0){\line(270,0){58}}
   \put(196,0){\line(270,0){83}}
   \put(45,-2.5){$\blacktriangleright$}
   \put(225,-2.5){$\blacktriangleright$}
   \put(133,23){$\blacktriangleright$}
   \put(133,-29){$\blacktriangleright$}
   \put(85,13){$\blacktriangleright$}
   \put(175,13){$\blacktriangleright$}
   \put(90,0){\circle*{3}}
   \put(180,0){\circle*{3}}
   \put(90,0){\circle{30}}
   \put(180,0){\circle{30}}
   \put(90,-10){$a$}
   \put(180,-10){$b$}
   \put(80,0){$\tilde U$}
   \put(170,0){$U$}
   \put(135, 30){$\mathnormal{\Sigma}_{1}$}
   \put(135, -40){$\mathnormal{\Sigma}_{2}$}
   \end{picture}
   \end{center}
   \caption{Contour $\Sigma_R$ for the Riemann-Hilbert problem for $R$}\label{6}
   \end{figure}

Using the parametrices $N$, $P$, and $\tilde{P}$, we
define the third transformation $S \mapsto R$ as follows
\begin{equation} \label{Rdef}
   R(z)= \left\{ \begin{array}{cl}
    S(z)N(z)^{-1} & \textrm{for }
    z \in \mathbb{C}\setminus\overline{(U\cup\tilde U)} \\
    S(z) P(z)^{-1} & \textrm{for } z \in U \\
    S(z)\tilde P(z)^{-1} & \textrm{for } z \in \tilde U
    \end{array} \right.
\end{equation}
Then $R$ has no jump on $[a,b]\setminus\overline{(U\cup\tilde U)}$, as
the jumps of $S$ and $N^{-1}$ cancel out. In $U$ and
$\tilde U$ the jumps of $S$ cancel out with the jumps of $P$ and
$\tilde P$, leaving only jumps for $R$ on the contour $\Sigma_R$ shown in Figure~\ref{6}.

The Riemann-Hilbert problem for $R$ is
\begin{align*}
\left\{\begin{array}{l} R(z)\textrm{ is analytic on }
\mathbb{C}\setminus\Sigma_R \\
R_{+}(z)=R_{-}(z)J_{R}(z)\textrm{ for } z\in\Sigma_R \\
R(z)=I+\mathcal{O}\left(\frac{1}{z}\right)\textrm{ for }
z\to\infty
\end{array}\right.
\end{align*}
where
\begin{align*}
J_{R}(z)=\left\{\begin{array}{l} N(z)J_{S}(z)N(z)^{-1}\textrm{ for }
z \in \Sigma_R \setminus (\partial U\cup\partial\tilde U) \\
P(z)N(z)^{-1}\textrm{ for } z \in \partial U \\
\tilde P(z)N(z)^{-1}\textrm{ for } z \in \partial\tilde U
\end{array}\right.
\end{align*}

The jump matrices $J_{R}(z)=N(z)J_{S}(z)N(z)^{-1}$ tend to the identity
matrix at an exponential rate as $n \to \infty$.
The jump matrices on $\partial U$ and $\partial \tilde U$ tend to
the identity matrix but at a slower rate of $1/n$ as $n \to \infty$.
The precise form is obtained from the asymptotic expansion
of the Airy function as $z \to \infty$, $-\pi < \arg z < \pi$,
(see \cite{Ho})
\begin{align}
\textrm{Ai}(z)\sim\frac{e^{-\frac{2}{3}z^{\frac{3}{2}}}}{2 \sqrt{\pi}
z^{\frac{1}{4}}}\sum\limits_{k=0}^{\infty}
\frac{(-1)^{k} \Gamma\left(3k+\frac{1}{2}\right)}{9^{k}(2k)! \Gamma\left(\frac{1}{2}\right)}
\frac{1}{z^{\frac{3}{2}k}}
\end{align}
and the corresponding asymptotic expansion for $\textrm{Ai}'(z)$.
Using these facts in the parametrix $P$ we find an asymptotic expansion
for the jump of $R$ on $\partial U$
\begin{equation} J_R(z) = P(z)N(z)^{-1} \sim
I+\sum\limits_{k=1}^{\infty}\frac{1}{n^{k}}\Delta_{k}(z)
 \label{As2} \end{equation}
where
\begin{align}
\Delta_{k}(z)= &
   \frac{1}{\sqrt{\pi}}
   \left(\frac{\Gamma\left(3k+\frac{1}{2}\right)}{9^{k}(2k)!}-
   \frac{\Gamma\left(3k-\frac{3}{2}\right)}{4\cdot9^{k-1}(2(k-1))!}\right)
   \frac{1}{\left(\frac{3}{2}\phi(z)\right)^{k}}I
   \nonumber\\
   &-\frac{1}{4\sqrt{\pi}}\frac{\Gamma\left(3k-\frac{3}{2}\right)}{9^{k-1}(2(k-1))!}
   \frac{1}{\left(\frac{3}{2}\phi(z)\right)^{k}}\sigma_{2}
   \quad \textrm{ for } k \textrm{ even}\label{Link1}
\end{align}
and
\begin{align}
\Delta_{k}(z)= &
   -\frac{\beta(z)^{2}}{\left(\frac{3}{2}\phi(z)\right)^{k}}\frac{1}{2\sqrt{\pi}}
   \left(\frac{\Gamma\left(3k+\frac{1}{2}\right)}{9^{k}(2k)!}-
   \frac{\Gamma\left(3k-\frac{3}{2}\right)}{2\cdot9^{k-1}(2(k-1))!}\right)
   \left(\sigma_{3}+i\sigma_{1}\right)\nonumber \\
   &-\frac{\beta(z)^{-2}}{\left(\frac{3}{2}\phi(z)\right)^{k}}
   \frac{1}{2\sqrt{\pi}}\frac{\Gamma\left(3k+\frac{1}{2}\right)}{9^{k}(2k)!}
   \left(\sigma_{3}-i\sigma_{1}\right)
    \quad \textrm{ for }k\textrm{ odd}\label{Link2}
\end{align}
where
\begin{equation} \label{Pauli}
    \sigma_{1} = \begin{pmatrix} 0 & 1 \\ 1 & 0 \end{pmatrix},
    \quad \sigma_{2} = \begin{pmatrix} 0 & -i \\ i & 0 \end{pmatrix},
    \quad \sigma_3 = \begin{pmatrix} 1 & 0 \\ 0 & -1 \end{pmatrix}
    \end{equation}
are the Pauli matrices.

A similar expansion
\begin{equation} J_R(z) = \tilde{P}(z) N(z)^{-1} \sim
    I + \sum_{k=1}^{\infty} \frac{1}{n^k} \tilde{\Delta}_k(z)
    \label{As3} \end{equation}
holds for the jump matrix on $\partial \tilde U$.

As a result we find by the methods of \cite{DKMVZ2},
see also \cite[Lemma 8.3]{KMVV},
\begin{lemma}
There exist matrix valued functions $R_k(z)$ with the
property that for every $l \in \mathbb N$, there exist constants $C > 0$ and $r > 0$
such that for every $z$ with $|z|\geq r$,
\begin{align}
\left\|R(z)-I-\sum\limits_{k=1}^{l}\frac{R_{k}(z)}{n^{k}}\right\|\leq\frac{C}{|z|n^{l+1}}
\end{align}
\end{lemma}

So we write
\begin{align}
R(z) \sim I+\sum\limits_{k=1}^{\infty}\frac{1}{n^{k}}R_{k}(z) \label{As1}
\end{align}
From (\ref{As1}), (\ref{As2}) and (\ref{As3}) and the Riemann-Hilbert
problem for $R$, we find an additive Riemann-Hilbert problem for $R_k(z)$,
\begin{align}\label{iterative}
\left\{\begin{array}{l} R_{k}(z) \textrm{ is analytic on }
    \mathbb{C}\setminus(\partial U \cup \partial \tilde U) \\
R_{k+}(z)=R_{k-}(z) +\sum\limits_{l=0}^{k-1}R_{l-}(z)\Delta_{k-l}(z)
    \textrm{ for } z \in \partial U \\
R_{k+}(z)=R_{k-}(z)+\sum\limits_{l=0}^{k-1}R_{l-}(z)\tilde\Delta_{k-l}(z)
    \textrm{ for } z \in \partial\tilde U \\
R_{k}(z)=\mathcal{O}\left(\frac{1}{z}\right)\textrm{ as } z\to\infty
\end{array}\right.
\end{align}
where $R_{0}(z)=I$. These Riemann-Hilbert problems can be
successively solved using the Sokhotskii-Plemelj formula, or using a
technique based on Laurent series expansion  as in
\cite{KMVV}.

\section{Proof of Theorem \ref{Main Theorem}}

For the proof of (\ref{annbnn}) we do not need to compute the
explicit forms of the $R_k$'s. However, we need to know that they
have the following structure. Recall that the Pauli matrices are
given in (\ref{Pauli}).

\begin{lemma}\label{LemmaR}
For $k$ odd, $R_{k}(z)$ is a linear combination of $\sigma_{1}$ and
$\sigma_{3}$ and for $k$ even, $R_{k}(z)$ is a linear combination of
$I$ and $\sigma_{2}$.
\end{lemma}
\begin{proof}
For $k=1$, we know because of (\ref{iterative}) that
$R_{1+}=R_{1-}+\Delta_{1}$ on $\partial U$ and
$R_{1+}=R_{1-}+\tilde\Delta_{1}$ on $\tilde\partial U$. As
$\Delta_{1}$,
$\tilde\Delta_{1}\in\textrm{span}\left\{\sigma_{1},\sigma_{3}\right\}$
on account of (\ref{Link2}), $R_{1}(z)$ must be a linear combination
of $\sigma_{1}$ and $\sigma_{3}$ as well.

Let $k\geq 1$ and once more observe (\ref{iterative}).
If $k$ is odd, then again by
(\ref{Link2}) $\Delta_{k}$,
$\tilde\Delta_{k}\in\textrm{span}\left\{\sigma_{1},\sigma_{3}\right\}$
and using induction on $k$, for  every $l < k$, $R_{l-}(z)\Delta_{k-l}(z)$
and $R_{l-}(z) \tilde\Delta_{k-l}(z)$ are products
of a linear combination of $\sigma_{1}$ and $\sigma_{3}$ and a
linear combination of $I$ and $\sigma_{2}$ (see also (\ref{Link1})--(\ref{Link2})),
which results in a
linear combination of $\sigma_{1}$ and $\sigma_{3}$.
Thus all terms in the (additive) jump for $R_k$ on $\partial U$ and on $\partial \tilde{U}$
are in the span of $\sigma_1$ and $\sigma_3$, and it follows that
$R_{k}\in\textrm{span}\left\{\sigma_{1},\sigma_{3}\right\}$ if $k$ is odd.

If $k$ is even, then by induction, where we use again
((\ref{Link1})--(\ref{Link2}), we have that $R_{l-}(z)
\Delta_{k-l}(z)$ and $R_{l-}(z) \tilde \Delta_{k-l}(z)$ are either
products of two linear combinations of $I$ and $\sigma_{2}$ (in case
$l$ is even), or products of two linear combinations of $\sigma_1$
and $\sigma_3$ (in case $l$ is odd). In both cases we find that
$R_{l-}(z) \Delta_{k-l}(z)$ and $R_{l-}(z) \tilde \Delta_{k-l}(z)$
are linear combinations of $I$ and $\sigma_2$, which implies that
$R_{k}\in\textrm{span}\left\{I,\sigma_{2}\right\}$ if $k$ is even.
\end{proof}

Now we can finally prove our main result.
\begin{proof}[Proof of Theorem \ref{Main Theorem}]
We start from the expressions (\ref{anNinY}) and (\ref{bnNinY}) for
$a_{n,n}$ and $b_{n,n}$ in terms of the solution of the
Riemann-Hilbert problem for $Y$. Following the transformations $Y
\mapsto T \mapsto S$, we find that
\begin{align} \label{anninS}
a_{n,n}=\left(S_{1}\right)_{12} \left(S_{1}\right)_{21}
\end{align}
and
\begin{align} \label{bnninS}
b_{n,n}=\frac{\left(S_{2}\right)_{12}}{\left(S_{1}\right)_{12}}-\left(S_{1}\right)_{22}
\end{align}
where $S_1$ and $S_2$ are the terms in the expansion of $S(z)$ as $z \to \infty$,
\[ S(z) = I + \frac{1}{z} S_1 + \frac{1}{z^2} S_2 + \mathcal{O}\left(\frac{1}{z^3}\right). \]
To obtain (\ref{bnninS}) we use that $g(z) = \log z + \mathcal{O}(1/z)$,
see also \cite{DK}.

By (\ref{Rdef}), we know that $S(z)=R(z)N(z)$ for $|z|$ large enough, so
we need the first terms in the expansions of $N(z)$ and $R(z)$ as $z \to \infty$.
From (\ref{Ndef}) we have
\begin{align}
N(z)&=\frac{\beta(z)+\beta(z)^{-1}}{2}I+\frac{\beta(z)-\beta(z)^{-1}}{2}\sigma_{2}\nonumber
\\
&=I-\frac{(b-a)}{4}\sigma_{2}\frac{1}{z}+
\left(\frac{(b-a)^{2}}{32}I-\frac{b^{2}-a^{2}}{8}\sigma_{2}\right)
\frac{1}{z^{2}}+\mathcal{O}\left(\frac{1}{z^{3}}\right) \label{Nexpansion}
\end{align}
and from Lemma~\ref{LemmaR}
\begin{align} \nonumber
R(z)&=I+ \frac{1}{z}\left(\sum\limits_{m\textrm{
odd}}\frac{1}{n^{m}}
\left(R_{m1\sigma_{1}}\sigma_{1}+R_{m1\sigma_{3}}\sigma_{3}\right)+\sum\limits_{m\textrm{
even}}\frac{1}{n^{m}}\left(R_{m1I}I+R_{m1\sigma_{2}}\sigma_{2}\right)\right)
\\
&+\frac{1}{z^{2}}\left(\sum\limits_{m\textrm{ odd}}\frac{1}{n^{m}}
\left(R_{m2\sigma_{1}}\sigma_{1}+R_{m2\sigma_{3}}\sigma_{3}\right)
+\sum\limits_{m\textrm{ even}}\frac{1}{n^{m}}
\left(R_{m2I}I+R_{m2\sigma_{2}}\sigma_{2}\right)\right)
\nonumber\\
&+\mathcal{O}\left(\frac{1}{z^{3}}\right) \label{Rexpansion}
\end{align}
where the constants $R_{mjI}$,  $R_{mj\sigma_{k}}$,
for $m\in \mathbb N$, $j=1,2$, and $k=1,2,3$ are such that
$R_{mjI} I + \sum_{k=1}^3 R_{mj\sigma_k} \sigma_k$ is the
coefficient of $z^{-j}$ in the Laurent expansion of $R_m(z)$
around $z=\infty$.

Therefore, by (\ref{Nexpansion}) and (\ref{Rexpansion}),
\begin{align}
S(z) = & R(z)N(z)\sim
I+\frac{1}{z}\left(-\frac{(b-a)}{4}\sigma_{2}+\sum\limits_{m\textrm{
odd}}\frac{1}{n^{m}}
\left(R_{m1\sigma_{1}}\sigma_{1}+R_{m1\sigma_{3}}\sigma_{3}\right)\right.
\nonumber\\
&\left.+\sum\limits_{m\textrm{
even}}\frac{1}{n^{m}}\left(R_{m1I}I+R_{m1\sigma_{2}}\sigma_{2}\right)\right)\nonumber\\
&+\frac{1}{z^{2}}\left(\frac{(b-a)^{2}}{32}I-\frac{b^{2}-a^{2}}{8}\sigma_{2}+
\sum\limits_{m\textrm{ odd}}\frac{1}{n^{m}}
\left(\left(R_{m2\sigma_{1}}+i\frac{b-a}{4}R_{m1\sigma_{3}}\right)\sigma_{1}\right.\right.\nonumber\\
&+\left.\left.\left(R_{m2\sigma_{3}}-i\frac{b-a}{4}R_{m1\sigma_{1}}\right)\sigma_{3}\right)+\sum\limits_{m\textrm{
even}}\frac{1}{n^{m}}
\left(\left(R_{m2I}-\frac{b-a}{4}R_{m1\sigma_{2}}\right)I\right.\right.\nonumber\\
&+\left.\left.\left(R_{m2\sigma_{2}}-\frac{b-a}{4}R_{m1I}\right)\sigma_{2}\right)\right)
+\mathcal{O}\left(\frac{1}{z^{3}}\right) \label{Sexpansion}
\end{align}
which implies that
\begin{align}\label{S12}
\left(S_{1}\right)_{12}\sim\frac{b-a}{4}i+\sum\limits_{m\textrm{
odd}}\frac{1}{n^{m}}R_{m1\sigma_{1}}-i\sum\limits_{m\textrm{
even}}\frac{1}{n^{m}}R_{m1\sigma_{2}}
\end{align}
and
\begin{align}\label{S21}
\left(S_{1}\right)_{21}\sim-\frac{b-a}{4}i+\sum\limits_{m\textrm{
odd}}\frac{1}{n^{m}}R_{m1\sigma_{1}}+i\sum\limits_{m\textrm{
even}}\frac{1}{n^{m}}R_{m1\sigma_{2}}\
\end{align}
Inserting (\ref{S12}) and (\ref{S21}) into (\ref{anninS}) then finally
gives
\begin{align*}
a_{n,n}\sim\frac{(b-a)^{2}}{16}+\sum\limits_{m=1}^{\infty} \frac{\alpha_{2m}}{n^{2m}}
\end{align*}
for certain constants $\alpha_{2m}$.

Similar to (\ref{S12}) and (\ref{S21}) we have that $(S_2)_{12}$ and $(S_1)_{22}$
have asymptotic expansions in powers of $1/n$.
From the expansion (\ref{Sexpansion}) for $S$, we see
\begin{align*}
\left(S_{2}\right)_{12}\sim & \frac{b^{2}-a^{2}}{8}i+\sum\limits_{m\textrm{
odd}}\frac{1}{n^{m}}\left(\frac{b-a}{4}iR_{m1\sigma_{3}}+R_{m2\sigma_{1}}\right) \nonumber\\
&+\sum\limits_{m\textrm{
even}}\frac{1}{n^{m}}i\left(\frac{b-a}{4}R_{m1I}-R_{m2\sigma_{2}}\right)
\end{align*}
and
\begin{align*}
\left(S_{1}\right)_{22}\sim-\sum\limits_{m\textrm{
odd}}\frac{1}{n^{m}}R_{m1\sigma_{3}}+\sum\limits_{m\textrm{
even}}\frac{1}{n^{m}}R_{m1I}
\end{align*}
From (\ref{bnninS}) it then follows that
\begin{align} \label{bnformule}
b_{n,n} \sim \sum_{m=0}^{\infty} \frac{\beta_m}{n^m}
\end{align}
where $\beta_0 = \frac{b+a}{2}$ and
\begin{align} \label{beta1formule}
    \beta_1 =
    2R_{11\sigma_{3}}-\frac{4}{b-a}iR_{12\sigma_{1}}+\frac{2(b+a)}{b-a}iR_{11\sigma_{1}}.
\end{align}

Our final task is to further evaluate the right-hand side of (\ref{beta1formule}).
As in \cite{KMVV}, we have that $\Delta_1$ is meromorphic in a neighborhood of
$b$ with a pole in $b$. Indeed, if we write
\begin{align} \label{B0def}
\frac{\beta(z)^{-2}}{\phi(z)}=(z-b)^{-2}\sum\limits_{m=0}^{\infty}B_{m}(z-b)^{m},
    \qquad B_0 = \frac{3}{2\pi h(b)},
\end{align}
and use (\ref{Link2}), then we find for $z$ in a neighborhood of $b$,
\begin{align} \nonumber
    \Delta_{1}(z)= &
     \left( - \frac{5B_1}{144 } \left(\sigma_3 - i \sigma_1\right)
    + \frac{7B_0}{144(b-a)} \left(\sigma_3 + i \sigma_1\right)\right) \frac{1}{z-b} \\
    &
    - \frac{5B_0}{144 } \left(\sigma_3 - i \sigma_1\right) \frac{1}{(z-b)^2}
    +\mathcal{O}\left(1\right).
    \label{Delta1expansion}
\end{align}
Similarly, for $z$ in a neighborhood of $a$, we have
\begin{align} \label{A0def}
    \frac{\beta(z)^2}{\tilde{\phi}(z)} = (z-a)^{-2} \sum_{m=0}^{\infty} A_m(z-a)^m,
    \qquad A_0 = \frac{3}{2\pi h(a)},
\end{align}
and
\begin{align} \nonumber
    \tilde{\Delta}_{1}(z)= &
    \left( - \frac{5A_1}{144 } \left(\sigma_3 + i \sigma_1\right)
    - \frac{7A_0}{144(b-a)} \left(\sigma_3 - i \sigma_1\right)\right) \frac{1}{z-a} \\
   &  - \frac{5A_0}{144 } \left(\sigma_3 + i \sigma_1\right) \frac{1}{(z-a)^2}
    +\mathcal{O}\left(1\right).
    \label{Deltatilde1expansion}
\end{align}
As in \cite{KMVV} we have that $R_1(z)$ for $z \in \mathbb C \setminus
\overline{U \cup \tilde{U}}$ is equal to the sum of the
Laurent parts of (\ref{Delta1expansion}) and (\ref{Deltatilde1expansion}).
Expanding $R_1(z)$ as $z \to \infty$, we then get
\begin{align*}
R_1(z) = R_{11} \frac{1}{z} + R_{12} \frac{1}{z^2} + \mathcal O \left(\frac{1}{z^3}\right)
\quad \textrm{as } z \to \infty,
\end{align*}
where
\begin{align*}
R_{11} = &
     - \frac{5A_1}{144 } \left(\sigma_3 + i \sigma_1\right)
    - \frac{7A_0}{144(b-a)} \left(\sigma_3 - i \sigma_1\right) \\
    &   - \frac{5B_1}{144 } \left(\sigma_3 - i \sigma_1\right)
    + \frac{7B_0}{144(b-a)} \left(\sigma_3 + i \sigma_1\right) \\
R_{12} = &  - \frac{5a A_1}{144 } \left(\sigma_3 + i \sigma_1\right)
    - \frac{7a A_0}{144(b-a)} \left(\sigma_3 - i \sigma_1\right)\\
    &  - \frac{5b B_1}{144 } \left(\sigma_3 - i \sigma_1\right)
    + \frac{7b B_0}{144(b-a)} \left(\sigma_3 + i \sigma_1\right) \\
    &  - \frac{5A_0}{144 } \left(\sigma_3 + i \sigma_1\right)
    - \frac{5B_0}{144 } \left(\sigma_3 - i \sigma_1\right).
\end{align*}
Thus
\begin{align} \label{R11sigma3}
    R_{11\sigma_3} & = - \frac{5(A_1 + B_1)}{144 } - \frac{7(A_0-B_0)}{144(b-a)}, \\
    R_{11\sigma_1} & = \label{R11sigma1}
     -i \frac{5(A_1-B_1)}{144} +i \frac{7(A_0+B_0)}{144(b-a)}, \\
    R_{12\sigma_1} & = - i \frac{5(a A_1-bB_1)}{144 } + i \frac{7(aA_0+bB_0)}{144(b-a)}
    - i \frac{5(A_0-B_0)}{144 }. \label{R12sigma1}
    \end{align}
Inserting (\ref{R11sigma3})--(\ref{R12sigma1}) into
(\ref{beta1formule}), we find after straightforward calculations
that $A_1$ and $B_1$ fully disappear  and that (\ref{beta1formule}) reduces to
\[ \beta_1 = \frac{B_0 - A_0}{3(b-a)}. \]
Using the explicit formulas for $A_0$ and $B_0$ given in (\ref{B0def})
and (\ref{A0def}), we
arrive at (\ref{beta1}), which completes the proof of Theorem \ref{Main Theorem}.
\end{proof}

\section*{Acknowledgement}
We are grateful to Pavel Bleher and Alexander Its for valuable
discussions on  their paper \cite{BI}.

\end{document}